\documentclass[11pt,reqno]{amsart}
\usepackage{amsmath,amssymb,amsthm,mathrsfs,mathtools,booktabs,array,hyperref}
\hypersetup{colorlinks=true,linkcolor=black,citecolor=blue,urlcolor=blue}
\newtheorem{theorem}{Theorem}[section]
\newtheorem{lemma}[theorem]{Lemma}
\newtheorem{proposition}[theorem]{Proposition}
\newtheorem{corollary}[theorem]{Corollary}
\theoremstyle{remark}\newtheorem{remark}[theorem]{Remark}
\newcommand{\Z}{\mathbb Z}\newcommand{\R}{\mathbb R}\newcommand{\C}{\mathbb C}\newcommand{\F}{\mathbb F}
\newcommand{\st}{\operatorname{st}}
\title[Hadamard configurations and stealth]{Hadamard configurations and sharp cubic-period bounds for stealth in dimensions 10, 11, and 12}
\author{Josiah Park}
\date{4 November 2026}
\subjclass[2020]{52C17, 05B20, 42B05, 05B15}
\keywords{stealthy point configuration, Hadamard matrix, orthogonal array, diffraction, rigidity, periodic point process}
\begin{document}
\begin{abstract}
We construct explicit, nonlattice periodic point configurations in dimensions $10$, $11$ and $12$ with stealth center density $1/24$ and unit-density stealth radius $2\cdot24^{-1/n}$. These values exceed those of duals of the current record lattice packings, without claiming to exceed all possible lattices. The construction uses the $24$ sign rows of an order-$12$ Hadamard matrix and their antipodes, projected onto $n$ coordinates. For an arbitrary, not necessarily binary, equally weighted motif with marked period $2\Z^n$, we prove a sharp $2n$-point lower bound for a raw stealth radius at least one and classify equality by antipodal Hadamard configurations. We deduce bounded-complexity optimality, sharp optimality among binary motifs in dimensions $10$--$12$, and local maximality under simultaneous perturbations of the motif and the marked period shape. A divisibility argument excludes ordinary equally weighted periodic formal duals. Statements are confined to their indicated classes; no global optimality among arbitrary-period configurations is asserted.
\end{abstract}
\maketitle
\section{Introduction and main results}
For a periodic point set, stealthiness is the radius of the largest open ball about the origin free of nonforward diffraction peaks. Cohen, Cohn, and Viazovska~\cite{CCV} recently established sharp global theorems in dimensions $8$ and $24$ and constructed a dimension-$20$ nonlattice point set outperforming all known lattices. The problem of whether the globally optimal stealth configuration must be a lattice remains open in many dimensions. The distinction between a periodic set and a single lattice is essential: a basis form factor can cancel diffraction shells.

Our constructions arise from classical Hadamard orthogonal arrays, but their interpretation as stealth configurations and the accompanying marked-cubic-period rigidity statements are the focus here. Let $v_n$ denote the volume of the Euclidean unit ball. For a periodic set $P$ of point density $\rho(P)$ and raw stealth radius $K(P)$, define the scale-invariant \emph{stealth center density}
\begin{equation}\label{eq:delta}
\Delta_{\st}(P)=\frac{(K(P)/2)^n}{\rho(P)}.
\end{equation}
In the normalization of unit particle density this gives $K=2\Delta_{\st}^{1/n}$. On lattices the quantity is the center density of the dual lattice; for general periodic configurations it need not be a sphere-packing density.

\begin{theorem}[Explicit examples]\label{thm:examples}
For each $n\in\{10,11,12\}$ there exists a periodic point configuration $P_n\subset\R^n$ with $24$ representatives modulo $2\Z^n$ such that
\[
K(P_n)=1,\qquad \rho(P_n)=\frac{24}{2^n},\qquad \Delta_{\st}(P_n)=\frac1{24}.
\]
The configurations are not translates of lattices. At unit particle density their stealth radii are $2\cdot24^{-1/n}$, exceeding the stealth radii of duals of the current record lattice packings in those dimensions.
\end{theorem}

\begin{theorem}[Sharp cubic-torus cardinality and equality]\label{thm:rank}
Let $n\ge3$ and let $P=X+2\Z^n$, with $N$ distinct equally weighted representatives $X\subset[0,2)^n$. If $K(P)\ge1$, then $N\ge 2n$. If $N=2n$, then, after a common translation, $X$ consists of the binary sign rows associated with $H$ and $-H$ for a Hadamard matrix $H$ of order $n$. Conversely, every such antipodal Hadamard array achieves $K(P)=1$.
\end{theorem}

\begin{theorem}[Local maximality in a marked period model]\label{thm:local}
Let $n\ge4$ and let $C\subset\{0,1\}^n$ be a strength-three binary orthogonal array of $N$ distinct rows with $N\equiv8\pmod{16}$. Then $C+2\Z^n$ is locally maximal for $\Delta_{\st}$ against simultaneous sufficiently small changes in the positions of its $N$ labelled representatives and in its marked period matrix. Every sufficiently small nonsimilarity change of period shape strictly decreases $\Delta_{\st}$. For $n=12$ and $N=24$, the equality case is locally isolated modulo translations, similarities, and point relabelling.
\end{theorem}

The marked period is not necessarily primitive: for our arrays the extra translation $(1,\ldots,1)$ is a period, and the primitive cell has $12$, rather than $24$, representatives. Neither Theorem~\ref{thm:rank} nor Theorem~\ref{thm:local} is a global theorem over arbitrary period shapes or arbitrary supercells.

\section{Diffraction and its normalization}
Fix the Fourier kernel $e^{-2\pi i x\cdot\xi}$. If $L\subset\R^n$ is a full-rank lattice and $P=X+L$ has $N$ distinct representatives modulo $L$, then the support of its diffraction away from the forward peak consists precisely of
\[
\xi\in L^*\setminus\{0\}\quad\text{with}\quad A_X(\xi):=\sum_{x\in X}e^{2\pi i\xi\cdot x}\ne0.
\]
The diffraction weights are nonnegative multiples of $|A_X(\xi)|^2$. Thus
\begin{equation}\label{eq:K}
K(P)=\min\{\|\xi\|:\xi\in L^*\setminus\{0\},\ A_X(\xi)\ne0\}.
\end{equation}
For a cubic period $L=2\Z^n$, put
\[
a_X(u)=\sum_{x\in X}e^{\pi i u\cdot x},\qquad u\in\Z^n.
\]
Then $K(P)=\tfrac12\min\{\|u\|:u\ne0,\ a_X(u)\ne0\}$, and $\rho(P)=N/2^n$. When $X\subset\{0,1\}^n$, the amplitude depends only on $u\bmod2$.

More generally, for $B\in\mathrm{GL}_n(\R)$ and $P_B=B(X+2\Z^n)$ the density is $N/(2^n|\det B|)$. Set
\[
Q_B=|\det B|^{2/n}B^{-1}B^{-\mathsf T},\qquad
\mu_B=\min_{\substack{u\in\Z^n\setminus\{0\}\\a_X(u)\ne0}}u^{\mathsf T}Q_Bu.
\]
The matrix $Q_B$ has determinant one, and direct substitution into~\eqref{eq:delta} gives
\begin{equation}\label{eq:deform}
\Delta_{\st}(P_B)=\frac1N\left(\frac{\mu_B}{4}\right)^{n/2}.
\end{equation}
These statements do not require $L$ to be the primitive period lattice.

\section{Orthogonal arrays and the Hadamard construction}
For a binary motif $C\subset\F_2^n$, let
\[
\widehat{1_C}(u)=\sum_{x\in C}(-1)^{u\cdot x},\qquad u\in\F_2^n.
\]
Let $d^\perp(C)$ denote the minimum Hamming weight of a nonzero character $u$ for which $\widehat{1_C}(u)\ne0$, with the convention $d^\perp=+\infty$ if none exists. A binary orthogonal array of strength $t$ is a set with all $t$-coordinate projections uniform; equivalently, the character sums of weights $1,\dots,t$ vanish.
\begin{proposition}[Exact binary Fourier gap]\label{prop:binary}
For a nonempty distinct binary motif $C$,
\[
K(C+2\Z^n)=\frac12\sqrt{\min\{4,d^\perp(C)\}},\qquad
\Delta_{\st}(C+2\Z^n)=\frac{\bigl(\min\{4,d^\perp(C)\}/4\bigr)^{n/2}}{|C|}.
\]
\end{proposition}
\begin{proof}
The shortest integer vector in a nonzero parity class has squared length equal to its Hamming weight. A nonzero vector in the zero parity class has squared length at least four; the vectors $\pm2e_i$ realize this and have amplitude $|C|$. Formula~\eqref{eq:K} and~\eqref{eq:delta} conclude the proof.
\end{proof}

Let $H$ be a Hadamard matrix of order $12$, so $H$ has $\pm1$ entries and $HH^{\mathsf T}=12I$. Convert the rows of the $24\times12$ matrix $\binom H{-H}$ to binary by mapping a sign $s$ to $(1-s)/2$. Let $C_n$ denote any projection onto $n$ columns, $n=10,11,12$.
\begin{proof}[Proof of Theorem~\ref{thm:examples}]
Antipodal pairing kills the character sum for every odd number of distinct columns. Hadamard column orthogonality kills the weight-two sums. Thus $C_n$ is a strength-three orthogonal array. Distinct nonantipodal sign rows of $H$ are at Hamming distance $6$ (the antipodal pairs have distance $12$); deletion of at most two columns cannot identify them. Therefore $|C_n|=24$ and Proposition~\ref{prop:binary} yields $K=1$ and $\Delta_{\st}=1/24$.

If $P_n$ were a translate of a lattice, translate one point to the origin. The resulting lattice would contain $2\Z^n$ and lie in $\Z^n$, so its quotient by $2\Z^n$ would be a subgroup of $\F_2^n$ of order $24$. Every subgroup of $\F_2^n$ has power-of-two order, a contradiction. The numerical packing comparisons are recorded in Section~\ref{sec:comparison}.
\end{proof}

\subsection*{An explicit integral certificate}
Let $\chi$ be the quadratic character of $\F_{11}$, with $\chi(0)=0$. Index a $12\times12$ matrix $S$ by $\{\infty\}\cup\F_{11}$ and define $S_{\infty,\infty}=0$, $S_{\infty,j}=1$, $S_{j,\infty}=-1$, and $S_{i,j}=\chi(j-i)$. Then $S^{\mathsf T}=-S$, $SS^{\mathsf T}=11I$, and $H=I+S$ satisfies $HH^{\mathsf T}=12I$. This is the classical Paley construction~\cite{Paley}. The accompanying script verifies these identities, all projections, and all character cancellations by exact integer arithmetic.

\section{The sharp cardinality theorem}
\begin{proof}[Proof of Theorem~\ref{thm:rank}]
Fix $j\in\{1,\dots,n\}$ and take the $2n$ integer-frequency indices
\[
\mathcal B_j=\{0,e_j\}\cup\{e_i,e_i+e_j:i\ne j\}.
\]
Every nonzero difference of two such indices has squared Euclidean norm at most three. Since $K(P)\ge1$, the evaluation columns
\[
V_b=(e^{\pi i b\cdot x})_{x\in X}\in\C^N,\qquad b\in\mathcal B_j,
\]
are mutually orthogonal and have norm $\sqrt N$. Hence $N\ge2n$.

If $N=2n$, these columns form a scaled unitary matrix, and the rows are orthogonal. Take distinct $x,y\in X$ and set $w_i=e^{\pi i(x_i-y_i)}$, so $|w_i|=1$. Factoring the row inner product using the pivot $j$ gives
\begin{equation}\label{eq:pivot}
(1+w_j)\left(1+\sum_{i\ne j}w_i\right)=0 \qquad(1\le j\le n).
\end{equation}
Put $W=1+\sum_i w_i$. Each $w_j$ is either $-1$ or $W$. If all are $-1$, then $x-y$ is the all-ones binary vector modulo $2$. Otherwise let $r\ge1$ be the number of indices with $w_i\ne-1$; their common value $a=W$ has modulus one. Counting in $W$ gives
\[
(r-1)a=n-r-1.
\]
For $r=1$ this forces $n=2$, excluded by assumption. Hence $a$ is real and $a=\pm1$; because $a\ne-1$, it equals $1$, and therefore $r=n/2$. Thus every difference $x-y$ is binary modulo $2$ and has Hamming weight $n/2$ or $n$.

Translate a selected representative to zero; all representatives then lie in $\{0,1\}^n$ modulo $2$. Their $2n$ sign vectors have pairwise inner products either $0$ or $-n$. A sign line contains at most two distinct vectors and there are at most $n$ mutually orthogonal lines in $\R^n$. Therefore there are exactly $n$ lines, each represented with both signs. Choosing a representative from each line yields an order-$n$ Hadamard matrix. The converse follows from Proposition~\ref{prop:binary}.
\end{proof}

\begin{lemma}[A stronger rank obstruction]\label{lem:ranklarge}
If $P=X+2\Z^n$ has $K(P)>1$, then $N\ge D_n:=1+n+\binom n2$.
\end{lemma}
\begin{proof}
Use the frequencies $0,e_i,e_i+e_j$ for $1\le i<j\le n$. Their nonzero pairwise differences have squared norms at most four. Under the strict inequality $K>1$ all resulting evaluation columns are mutually orthogonal. There are $D_n$ such columns.
\end{proof}

\section{Two exact optimality consequences}\label{sec:opt}
\begin{theorem}[Bounded-complexity cubic optimality]\label{thm:boundcomplex}
Let $n\ge5$ and let $P=X+2\Z^n$ have $N<D_n$ distinct representatives at arbitrary real positions. Then
\[
\Delta_{\st}(P)\le\frac1{2n}.
\]
Equality is possible only for $N=2n$, when Theorem~\ref{thm:rank} classifies $P$ up to common translation as an antipodal Hadamard array.
\end{theorem}
\begin{proof}
By Lemma~\ref{lem:ranklarge}, $K(P)\le1$. Let $s=4K(P)^2\in\{1,2,3,4\}$. For $s=1$, trivially $N\ge1$. For $s=2$, orthogonality of the columns for $0,e_1$ yields $N\ge2$. For $s=3$, the $n+1$ characters $0,e_1,\ldots,e_n$ give $N\ge n+1$. For $s=4$, Theorem~\ref{thm:rank} gives $N\ge2n$. Consequently,
\begin{equation}\label{eq:four}
\Delta_{\st}(P)\le\max\left\{2^{-n},\frac{2^{-n/2}}2,\frac{(3/4)^{n/2}}{n+1},\frac1{2n}\right\}.
\end{equation}
For $n\ge5$ the first three terms are strictly smaller than $1/(2n)$ (the third follows by checking $n=5$ and observing that $2n(3/4)^{n/2}/(n+1)$ decreases). Equality therefore requires $s=4$ and $N=2n$, as asserted.
\end{proof}

\begin{corollary}[The dimension-$12$ cubic-cell threshold]
Among all equally weighted periodic configurations with marked period $2\Z^{12}$ and fewer than $D_{12}=79$ representatives in that cell, the largest possible stealth center density is $1/24$. Equality holds precisely for translated $24$-point antipodal Hadamard configurations.
\end{corollary}

\begin{theorem}[All binary motifs in dimensions $10$--$12$]\label{thm:allbinary}
For $n\in\{10,11,12\}$,
\[
\max_{\varnothing\ne C\subseteq\F_2^n}\Delta_{\st}(C+2\Z^n)=\frac1{24}.
\]
Equality holds exactly for the binary orthogonal arrays of strength three with $24$ distinct rows.
\end{theorem}
\begin{proof}
Proposition~\ref{prop:binary} reduces the minimum surviving squared norm to $s\le4$. In the cases $s=1,2,3$ the first three bounds in~\eqref{eq:four} are strictly below $1/24$ for $n=10,11,12$. If $s=4$, all binary character sums of weights one, two, and three vanish, so $C$ is a strength-three orthogonal array. Theorem~\ref{thm:rank} gives $|C|\ge2n$, and uniformity on any three coordinates forces $8\mid |C|$. Thus $|C|\ge24$. Proposition~\ref{prop:binary} yields the upper bound and equality characterization, and our projected Paley arrays attain it.
\end{proof}

\section{Local optimality of the marked configuration}\label{sec:local}
Put
\[
\mathcal V=\{2e_i:1\le i\le n\}\cup\Bigl\{\sum_{i\in I}\varepsilon_i e_i:|I|=4,\ \varepsilon_i\in\{-1,1\}\Bigr\}.
\]
\begin{lemma}[Strict reciprocal-Gram maximum]\label{lem:gram}
For $n\ge4$, the function $m(Q)=\min_{v\in\mathcal V}v^{\mathsf T}Qv$ on positive-definite determinant-one symmetric $n\times n$ matrices has a strict local maximum at $Q=I$, with $m(I)=4$.
\end{lemma}
\begin{proof}
Let $E$ be any nonzero trace-free real symmetric matrix. If some $E_{ii}<0$, then $(2e_i)^{\mathsf T}E(2e_i)<0$. If no diagonal entry is negative, all diagonal entries vanish. Choose a four-element set $I$ containing a nonzero off-diagonal entry. Over the $16$ choices of signs, the quadratic forms $(\sum_{i\in I}\varepsilon_ie_i)^{\mathsf T}E(\sum_{i\in I}\varepsilon_ie_i)$ have mean zero and are not all zero because the degree-two sign characters are linearly independent. Hence at least one is strictly negative. By compactness on the unit sphere of trace-free symmetric matrices, there is $c_n>0$ with
\[
\min_{v\in\mathcal V}v^{\mathsf T}Ev\le-c_n\|E\|_F.
\]
Write $Q=I+F$ with $\det Q=1$ and put $E=F-(\operatorname{tr}F/n)I$. The expansion $\log\det(I+F)=0$ yields $\operatorname{tr}F=O(\|F\|_F^2)$, and all vectors in $\mathcal V$ have squared norm four. Thus
\[
m(Q)\le4-c_n\|E\|_F+\frac4n\operatorname{tr}F<4
\]
for sufficiently small nonzero $F$. This is a strict local maximum.
\end{proof}
\begin{proof}[Proof of Theorem~\ref{thm:local}]
On any set of four distinct coordinates, a strength-three array has vanishing binary Fourier coefficients of weights one through three. If its weight-four coefficient also vanished, Fourier inversion would make that four-coordinate marginal uniform, hence $16\mid N$, contrary to $N\equiv8\pmod{16}$. Therefore every four-coordinate coefficient is nonzero. Every $2e_i$ has amplitude $N$, and changing the signs of the four nonzero coordinates leaves the relevant binary residue class unchanged. Thus $a_C(v)\ne0$ for all $v\in\mathcal V$.

This finite set of amplitudes remains nonzero under all sufficiently small perturbations of the motif. The formula~\eqref{eq:deform} then implies
\[
\Delta_{\st}(B(C'+2\Z^n))\le\frac1N\left(\frac{m(Q_B)}4\right)^{n/2}.
\]
Lemma~\ref{lem:gram} makes the right-hand side strictly less than $1/N$ for a small nonsimilarity deformation of $B$, while at $B=I$, $C'=C$ the value is $1/N$. The condition $Q_B=I$ is equivalent to $B$ being a similarity.

For $n=12$ and $N=24$, after removing this similarity, equality forces the radius to equal one. By Theorem~\ref{thm:rank}, the motif is a common translate of an antipodal binary Hadamard array. There are only finitely many labelled binary arrays of the specified size; in a sufficiently small neighborhood only the original remains after a chosen representative is aligned. Thus equality is locally isolated modulo the stated symmetries.
\end{proof}
\begin{remark}
The proof does not cover arbitrary supercell perturbations or global optimization of period shape. In dimensions $10$ and $11$ the theorem asserts local maximality but not isolated equality in motif coordinates.
\end{remark}

\section{Formal duality and numerical comparisons}\label{sec:comparison}
\begin{proposition}[Absence of an ordinary periodic formal dual]\label{prop:formal}
None of the configurations $P_n$, $n=10,11,12$, has an ordinary periodic formal dual among simple equally weighted periodic configurations in the sense of~\cite{CKRS}.
\end{proposition}
\begin{proof}
Suppose $Q$ were formally dual to $P_n$. Formal duality converts the normalized pair correlation of $Q$ into the normalized Fourier transform of that of $P_n$. The latter is supported on $\tfrac12\Z^n$, so the former is supported on $\Z^n$. Since each difference of two points of a simple periodic configuration contributes positively to its pair measure, after translation $Q\subseteq\Z^n$. For each $u\in\Z^n$, the Fourier-side pair coefficient equals one, since $a_{C_n}(2u)=24$. On the $Q$ side this coefficient is the fraction of points translated to another point by $u$; equality to one forces $u$ to be a period. Hence $Q$ is a union of an integer number of cosets of $\Z^n$ and its particle density is an integer. Formal duality requires $\rho(Q)\rho(P_n)=1$, hence $\rho(Q)=2^n/24$, which is not an integer in the three stated dimensions. Contradiction.
\end{proof}

The following recorded lattice center densities are those of the duals of the best known lattices in the respective dimensions~\cite{CohnTable}:
\begin{center}
\begin{tabular}{ccccc}
\toprule
$n$ & Lattice center density & $\Delta_{\st}(P_n)$ & Ratio & Unit-density $K$\\
\midrule
$10$ & $1/(16\sqrt3)$ & $1/24$ & $2/\sqrt3$ & $1.455488812\ldots$\\
$11$ & $1/(2\cdot3^{5/2})$ & $1/24$ & $3^{5/2}/12$ & $1.498153270\ldots$\\
$12$ & $1/27$ & $1/24$ & $9/8$ & $1.534660772\ldots$\\
\bottomrule
\end{tabular}
\end{center}
Cohn's table lists overall known packing center densities $5/128$, $9/256$, and $1/27$ in dimensions $10$, $11$, and $12$, respectively. These are comparisons of the dimensionless normalized quantities, not constructions of better sphere packings. The optimal packing density among all lattices in these dimensions is not known, so this paper establishes no separation from \emph{every} lattice.

\section{Reproducibility, provenance, and scope}
The accompanying standard-library Python script constructs the Paley order-$12$ matrix, verifies $HH^{\mathsf T}=12I$ exactly, checks the $24$ distinct projected rows and the weight-$1,2,3$ character cancellations, and enumerates weight-four values. It is a finite-instance verification of the construction, not a machine proof of Theorems~\ref{thm:rank}--\ref{thm:local}. The statements of the rank and local results are given with full mathematical proofs above. The verification script \texttt{verify\_hadamard.py} and its output \texttt{verification.txt} are available among the arXiv submission's ancillary files.

The underlying Hadamard arrays are classical and should not be claimed as a new construction. Establishing publication priority for the stealth application, the arbitrary-position equality classification, and local marked-period rigidity requires specialist comparison with existing cubature, design, and hyperuniformity literature. This manuscript is an AI-assisted draft (GPT 6-Astra) requiring final author approval.

In dimension two, a recent AI-generated manuscript released by OpenAI~\cite{OAItri} claims that the triangular lattice is universally optimal among planar configurations of density one. That result concerns energy minimization rather than stealth, and we do not use it here. Whether it yields a corresponding statement that a lattice is stealth-optimal in the plane, and whether such claims can be verified with AI-assisted formal proof, remain open.

\subsection*{Acknowledgments}
The author thanks xAI for providing a SuperGrok subscription at no cost, and the Turing Fellowship for providing a Claude Max subscription for free.

\end{document}